\documentclass[a4paper,12pt,leqno]{article}
\usepackage{amsmath,amsfonts,amsthm,amssymb,dsfont}
\usepackage[alphabetic]{amsrefs}
\usepackage[OT4]{fontenc}
\usepackage{enumerate}
\usepackage{epsfig}
\usepackage{epstopdf}

\long\def\symbolfootnote[#1]#2{\begingroup
\def\thefootnote{\fnsymbol{footnote}}\footnote[#1]{#2}\endgroup}

\newtheorem{theorem}{Theorem}[section]
\newtheorem{lemma}[theorem]{Lemma}

\newtheorem{prop}[theorem]{Proposition}
\newtheorem{cor}[theorem]{Corollary}

\theoremstyle{definition}
\newtheorem{rem}[theorem]{Remark}

\newtheorem{defin}[theorem]{Definition}

\newcommand {\Z} {{\mathbb{Z}}}

\renewcommand {\S} {{\mathbb{S}}}
\newcommand {\Mn}{MS(L)}
\newcommand {\In}{IS(L)}
\newcommand {\MVn}{\mathcal{MS}(L)}

\newcommand {\M}{MS(E,\gamma,\alpha)}
\newcommand {\MV}{\mathcal{MS}(E,\gamma,\alpha)}
\newcommand {\I}{IS(E,\gamma,\alpha)}

\begin{document}

\begin{center}
\large\bfseries Contractibility of the Kakimizu complex and symmetric Seifert surfaces
\end{center}

\begin{center}\bf
Piotr Przytycki$^a$\symbolfootnote[1]{Partially supported by MNiSW grant
N201 012 32/0718, MNiSW grant
N N201 541738 and the Foundation for Polish Science.}  \& Jennifer Schultens$^b$\symbolfootnote[2]{Partially supported by NSF grant.}
\end{center}

\begin{center}\it
$^a$ Institute of Mathematics, Polish Academy of Sciences,\\
 \'Sniadeckich 8, 00-956 Warsaw, Poland\\
\emph{e-mail:}\texttt{pprzytyc@mimuw.edu.pl}
\end{center}

\begin{center}\it
$^b$ Department of Mathematics, One Shields Avenue,\\
University of California, Davis, CA 95616, USA\\
\emph{e-mail:}\texttt{jcs@math.ucdavis.edu}
\end{center}

\begin{abstract}
\noindent Kakimizu complex of a knot is a flag simplicial complex
whose vertices correspond to minimal genus Seifert surfaces and
edges to disjoint pairs of such surfaces. We discuss a general
setting in which one can define a similar complex. We prove that
this complex is contractible, which was conjectured by Kakimizu.
More generally, the fixed-point set (in the Kakimizu complex) for
any subgroup of an appropriate mapping class group is contractible
or empty. Moreover, we prove that this fixed-point set is
non-empty for finite subgroups, which implies the existence of
symmetric Seifert surfaces.
\end{abstract}

\section{Introduction}

We study a generalisation $MS(E)$ of the following simplicial
complex $\Mn$ defined by Kakimizu \cite{K}. Let $E=E(L)$ be the
exterior of a tubular neighbourhood of a knot $L$ in $\S^3$. A
\emph{spanning surface} is a surface properly embedded in $E$,
which is contained in some Seifert surface for $L$. Let $\MVn$ be
the set of isotopy classes of spanning surfaces which have minimal
genus. The vertex set of $\Mn$ is defined to be $\MVn$. Vertices
$\sigma, \sigma'\in \MVn$ span an edge if they have representative
spanning surfaces which are disjoint. Simplices are spanned on all
complete subgraphs of the $1$--skeleton. In other words, $\Mn$ is
the \emph{flag} complex spanned on its $1$--skeleton. Kakimizu
defines $\Mn$ for links in the same way, but we later argue that
this is not the right definition and we define our $MS(E)$ for
$E=E(L)$ differently. However, for all links whose $\Mn$ have been
so far studied we have $MS(E(L))=\Mn$.

The general setting in which we define $MS(E(L))$, or more
generally $\M$, is the following. Let $E$ be a compact connected
orientable, irreducible and $\partial$--irreducible $3$--manifold.
In particular, for any non-splittable link $L$ in $\S^3$, the
complement $E(L)$ of a regular neighbourhood of $L$ satisfies
these conditions. Let $\gamma$ be a union of oriented disjoint
simple closed curves on $\partial E$, which does not separate any
component of $\partial E$. For $E=E(L)$ an example of $\gamma$ is
the set of longitudes of all link components (or its subset). We
fix a class $\alpha$ in the homology group $H_2(E,\partial E,\Z)$
satisfying $\partial \alpha = [\gamma]$. For $E=E(L)$ and $\gamma$
the set of longitudes, there is only one choice for $\alpha$. It
is the homology class dual to the element of $H^1(E,\Z)$ mapping
all oriented meridian classes onto a fixed generator of $\Z$. A
\emph{spanning surface} is an oriented surface properly embedded
in $E$ in the homology class $\alpha$ whose boundary is homotopic
with $\gamma$.

We now define the simplicial complex $\M$, which we abbreviate to
$MS(E)$, if $E=E(L)$ and $\gamma$ is the set of all longitudes.
The vertex set of $\M$ is defined to be $\MV$, the set of isotopy
classes of spanning surfaces which have minimal genus. However, we
span an edge on $\sigma, \sigma'\in \MV$ only if they have
representatives $S\in\sigma,S'\in \sigma'$ such that the
(connected) lift of $E\setminus S'$ to the infinite cyclic cover
associated with $\alpha$ intersects exactly two lifts of
$E\setminus S$. In the terminology of Section~\ref{sec:distance}
this means that the \emph{Kakimizu distance} between $\sigma$ and
$\sigma'$ equals one. This is not always true for disjoint $S,S'$
(because they are allowed to be disconnected). This error is made
by Kakimizu \cite[formula 1.3(b)]{K} who does not distinguish
between $\Mn$ and $MS(E(L))$. However, both his and our article
prove that the right complex to consider is $MS(E(L))$.

For every link $L$ it is a basic question to determine the complex
$MS(E(L))$ which encodes the structure of the set of all minimal
genus spanning surfaces. This has been done for all prime knots of
at most 10 crossings by Kakimizu \cite[Theorem A]{K2}. Moreover,
questions about common properties of all $MS(E(L))$ (or rather
$\Mn$) have been asked. Here is a brief summary (for a broader
account, see \cite{Pe}).

Scharlemann--Thompson proved \cite[Proposition 5]{ST} that
$MS(E(L))$ is connected, in the case where $L$ is a knot. Later
Kakimizu \cite[Theorem A]{K} provided another proof for links.
Schultens \cite[Theorem 6]{S} proved that, in the case where $L$
is a knot, $MS(E(L))$ is simply-connected (see also \cite{SS} for
atoroidal genus $1$ knots). For atoroidal knots bounds on the
diameter of $MS(E(L))$ have been obtained (\cite{Pe,SS}). Kakimizu
conjectured (see \cite[Conjecture 0.2]{Sa}) that $\Mn$ is
contractible. This was verified for special arborescent links by
Sakuma \cite[{Theorem 3.3 and Proposition 3.11}]{Sa}, and
announced for special prime alternating links by Hirasawa--Sakuma
\cite{HS}. In the present article, we confirm this conjecture,
under no hypothesis, for the complex $\M$.

\begin{theorem}
\label{thm:Contractibility}
$\M$ is contractible.
\end{theorem}

Using the same method we are also able to establish the following.
Note that for $E=E(L)$ all mapping classes of $E$ fix $\alpha$ and
the homotopy class of $\gamma$.

\begin{theorem}
\label{thm:fixed point} Let $G$ be a finite subgroup of the
mapping class group of $E$ fixing $\alpha$ and the homotopy class
of $\gamma$. We consider its natural action on $\M$. Then there is
a simplex in $\M$ fixed by all elements of $G$.
\end{theorem}

Sakuma argued \cite[Proposition 4.9(1)]{Sa} (see also
\cite[Theorem 5]{S} for knots) that the set of vertices of any
simplex of $\M$ can be realised as a union of pairwise disjoint
spanning surfaces. Hence in the language of spanning surfaces
Theorem~\ref{thm:fixed point} amounts to the following.

\begin{cor}
\label{cor:invariant surfaces} Let $G$ be a finite subgroup of the
mapping class group of $E$ fixing $\alpha$ and the homotopy class
of $\gamma$. Then there is a union of pairwise disjoint spanning
surfaces of minimal genus which is $G$--invariant up to isotopy.
\end{cor}

In the case where $E$ is atoroidal and $\partial E$ is a union of
tori, its interior admits, by the work of Thurston and the
theorem of Prasad, a unique complete hyperbolic structure. Then
the mapping class group of $E$ coincides with the isometry group
of its interior, hence it is finite. Moreover, after deforming the
metric in a way discussed in \cite[Chapter 10]{Pe} we can assume
that each element of $\MV$ has a unique representative of minimal
area. In this case Corollary~\ref{cor:invariant surfaces} gives
the following.

\begin{cor}
\label{cor:one surface} If $E$ is atoroidal and $\partial E$ is a
union of tori, then there is a union of pairwise disjoint
spanning surfaces of minimal genus which is invariant under any
isometry fixing $\alpha$ (the homotopy class of $\gamma$ is then
fixed automatically). In particular, if $E=E(L)$, then this union
is invariant under any isometry.
\end{cor}

A~related result concerning periodic knots was proved in Edmonds
\cite{E}.

Finally, Theorem~\ref{thm:Contractibility} turns out to be a special case ($G$ trivial) of the following.

\begin{theorem}
\label{thm:all contractible} Let $G$ be any subgroup of the
mapping class group of $E$ fixing $\alpha$ and the homotopy class
of $\gamma$. Then its fixed-point set $\mathrm{Fix}_G(\M)$  is
either empty or contractible.
\end{theorem}

We decided to provide first the proof of
Theorem~\ref{thm:Contractibility} and then the more technically
involved proof of the generalisation, Theorem~\ref{thm:all
contractible}.

We conclude with the following consequence of Theorem~\ref{thm:all contractible}.

\begin{cor}
\label{cor:EG} Denote by $G$ the mapping class group of $E$ fixing
$\alpha$ and the homotopy class of $\gamma$. Let $\mathcal F$ be
the set of those subgroups of $G$ which stabilise a point in $\M$.
Then $\M$ is the model for $E_{\mathcal F}(G)$ (the classifying
space for $G$ with respect to the family $\mathcal F$, see
\cite{L}).
\end{cor}

Actually, it is not clear to us what groups, apart from all finite
ones (see Theorem~\ref{thm:fixed point}), belong to the family
$\mathcal F$. It is also not clear if $\M$ can be locally
infinite.

\medskip

\textbf{Outline of the idea.} We now outline the main idea of the
article. The central object is the \emph{projection map} $\pi$,
which assigns to a pair of vertices $\sigma, \rho\in\MV$ at
distance $d>0$ a vertex $\pi_\sigma(\rho)$ adjacent to $\rho$ at
distance $d-1$ from $\sigma$. Kakimizu \cite{K} used the
projection to prove that $MS(E(L))$ is connected, but in fact he
did not need to verify that it is well-defined --- he worked only
with representatives of vertices. We verify that $\pi$ is
well-defined using a result of Oertel on \emph{cut-and-paste
operations} on surfaces with \emph{simplified intersection}.

We explain how to prove contractibility of $\M$. Assume for
simplicity that $\M$ is finite (which is the case for $E$
hyperbolic, see \cite[Corollary 8.8.6(b)]{Th}). We fix some $\sigma\in
\MV$. Then we prove that among vertices farthest from $\sigma$
there exists a vertex $\rho$ which is \emph{strongly dominated} by
$\pi_\sigma(\rho)$. This means that all the vertices adjacent to
$\rho$ are also adjacent to or equal $\pi_\sigma(\rho)$. Hence
there is a homotopy retraction of $\M$ onto the subcomplex spanned
by all the vertices except $\rho$. Proceeding in this way we
retract the whole complex onto $\sigma$.
\medskip

\textbf{Remaining questions.} Finally, we indicate that questions
about the structure of the set of all incompressible spanning
surfaces remain open. Kakimizu \cite{K} considers the complex
$\In$ whose vertices are isotopy classes of spanning surfaces
which are incompressible and $\partial$--incompressible but not
necessarily of minimal genus. The edges of $\In$ are defined like
edges of $\Mn$, in particular we have an embedding of $\Mn$ into
$\In$. Kakimizu asks if $\In$ is contractible as well. He proves
that $\In$ is connected, using a composition of the projection
$\pi$ with an additional operation, which we do not know how to
make well-defined on the set of isotopy classes of surfaces. This
is why we do not know if we can extend Theorem~\ref{thm:all
contractible} or even Theorem~\ref{thm:Contractibility} to the
complex $\In$ (or rather to $\I$, appropriately defined). Note
however that, since $\M$ would be a subcomplex of $\I$,
Theorem~\ref{thm:fixed point} would trivially carry over to $\I$.

\medskip

\textbf{Organisation of the article.} In
Section~\ref{sec:distance} we discuss \emph{Kakimizu distance}, a
geometric way to understand the distance between vertices of $\M$
in its $1$--skeleton. In Section~\ref{sec:simplified} we prove
that we can compute this distance from representative surfaces
with \emph{simplified intersection}. We use that in
Section~\ref{sec:projection} to prove that the projection map is
well-defined. In Section~\ref{sec:order} we introduce the order on
$\MV$ in which we will contract the complex. We establish various
properties of the projection map in Section~\ref{sec:properties}.
Using these, we establish contractibility,
Theorem~\ref{thm:Contractibility}, in
Section~\ref{sec:contractibility}. Next, in
Section~\ref{sec:fixed} we prove the fixed-point result,
Theorem~\ref{thm:fixed point}. Finally, in Section~\ref{sec:eg} we
prove Theorem~\ref{thm:all contractible} that all fixed-point sets
are contractible, if non-empty.
\medskip

\textbf{Acknowledgements.} After having proved
Theorem~\ref{thm:Contractibility}, we learned that Victor Chepoi
has outlined independently a possibly similar proof. In fact, our
article is inspired by what we have learned from \cite{CO} and
\cite{Pol}. We were also inspired by an argument which we have
learned from Saul Schleimer proving contractibility of the arc
complex. We thank Saul Schleimer for advice, encouraging us to
prove Theorem~\ref{thm:fixed point} and for telling us about
\cite{Pe}. We thank Jessica Banks for pointing out an error in our 
previous definition of semi-convexity. The first author is grateful to the Hausdorff Institute
of Mathematics in Bonn and to the Erwin Schr\"odinger Institute in
Vienna. The second author is grateful to the Max-Planck Institute
in Bonn.

\section{Kakimizu distance}
\label{sec:distance} In this section we start recalling the method
using which Kakimizu proved \cite[Theorem A]{K} that $MS(E(L))$ is
connected. This method was later used by Schultens \cite[Theorem
6]{S} to prove that $MS(E(L))$ is simply connected, in the case
where $L$ is a knot, and will be also the basic tool in the
present article.

The method is to study a pair $S,R$ of spanning surfaces via the
lifts of $E\setminus S, E\setminus R$ to the infinite cyclic cover
$\widetilde{E}$ of $E$ associated with the (kernel of the) element
of $H^1(E,\Z)$ dual to $\alpha$. It turns out that the distance in
$\M$ between two vertices $[S],[R]$ determined by those surfaces
can be read instantly from the relative position of the lifts of
$E\setminus S$ and $E\setminus R$.
\medskip

We recall the setting and notation of \cite{K}. Let $p\colon \widetilde{E}\rightarrow E$ be the covering map discussed above. Let $\tau$ be the generator of the group of covering transformations of $\widetilde{E}$. Suppose that $S\subset E$ is a spanning surface.
The hypothesis that $\gamma$ does not separate the components of $\partial E$ guarantees that $E\setminus S$ is connected. Let $E_0$ denote a lift of $E\setminus S$ to $\widetilde{E}$ and denote $E_j=\tau^j(E_0)$ for $j\in \Z$. Note the difference with \cite{K}, where $E_0$ is the closure of our $E_0$. Denote also $S_j=\overline{E}_{j-1}\cap \overline{E}_j$ for $j\in \Z$ (the bars will always denote closures).

\begin{figure}[htbp]
\begin{center}
\includegraphics[width=4in]{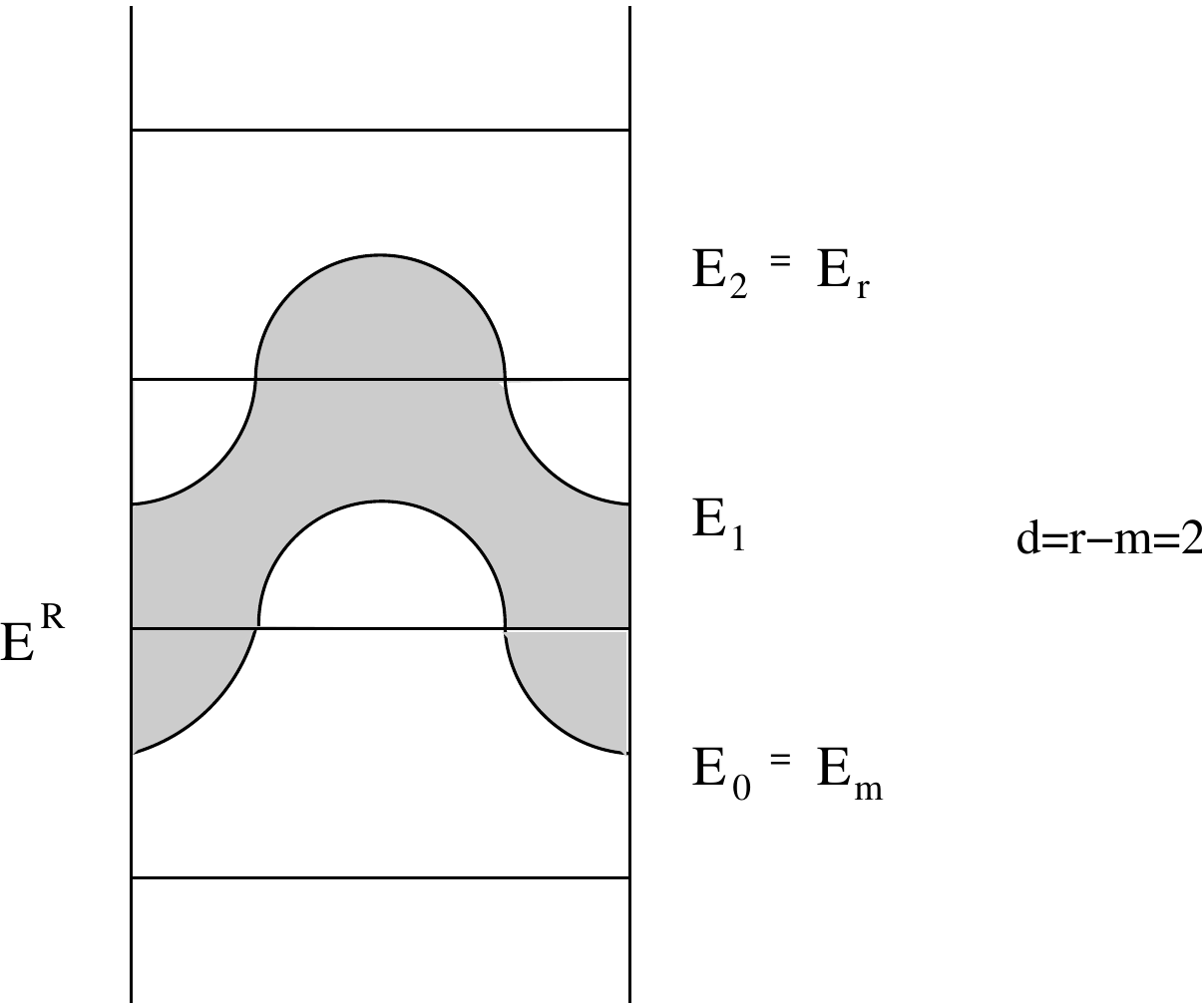}
\caption{$d(S,R)$ is defined via the lifts of $S$ and $R$}
\end{center}
\end{figure}

\begin{defin}
\label{def:distance} Let $R$ be another spanning surface. Let
$E^R$ be any lift of $E\setminus R$ to $\widetilde{E}$. We set
$$r=\max\{k\in \Z| E_k\text{ intersects } E^R\},\ m=\min\{k\in \Z|
E_k\text{ intersects } E^R\}$$ and we put $d(S,R)=r-m$. This value
does not depend on the choice of the lift $E^R$. See Figure 1.

Furthermore, for any two isotopy classes $\sigma, \rho$ of
spanning surfaces we define $d(\sigma, \rho)$ to be the minimum of
$d(S,R)$ over all representatives $S$ of $\sigma$ and $R$ of
$\rho$.
\end{defin}

Observe that in the case $\sigma=\rho$ we can take $S=R$ which satisfy $d(S,R)=0$. Recall that we declared two different vertices $\sigma,\rho$ of $\M$ to be adjacent if they satisfy $d(S,R)=1$ for some $S\in\sigma, R\in \rho$. Note that if $S$ and $R$ are disconnected, it could happen that $S$ and $R$ are disjoint, but $d(S,R)$ exceeds $1$. One might not be able to improve that by varying $S$ and $R$ in the isotopy classes.

Kakimizu proves the following. (Our context is more general, but the proof trivially carries over.)

\begin{prop}[{\cite[Proposition 1.4]{K}}]
\label{prop:distance}
The function $d$ is a metric on $\MV$.
\end{prop}

In fact, if we endow the $1$--skeleton of $\M$  with the path-metric $l$ in which all the edges have length $1$, then $d$ satisfies the following.

\begin{prop}[{\cite[Proposition 3.1]{K}}]
\label{prop:distances coincide}
The metric $d$ coincides with $l$ on $\MV$.
\end{prop}

Let us indicate how Kakimizu proves Proposition~\ref{prop:distances coincide}. The distance
$l=l(\sigma, \rho)$ is realised by a path $\sigma_0=\sigma, \sigma_1,\ldots, \sigma_l=\rho$.
By Proposition~\ref{prop:distance}, we have $d(\sigma,\rho)\leq d(\sigma_0,\sigma_1)+\ldots +d(\sigma_{l-1},\sigma_l)=l$, which is the estimate in one direction. The second estimate will be explained at the beginning of Section~\ref{sec:projection}.

\section{Simplified intersection}
\label{sec:simplified}
In this section we address the following issue. What hypotheses on the representatives $S,R$ of spanning surfaces $\sigma, \rho$ guarantee $d(\sigma, \rho)=d(S,R)$? To formulate a criterion we need the following terminology (see \cite{O}).

Let $S,R$ be compact surfaces properly embedded in a connected (not necessarily compact) $3$--manifold $M$ with boundary. We discuss \emph{product regions} bounded by $S$ and $R$ in $\partial M$ and $M$. If $\beta$ is an (abstract) arc, we denote by $\mathfrak I$ the product $\beta \times I$ with $\{x\} \times I$ collapsed to a point for each $x\in \partial \beta$. A \emph{product region} in $\partial M$ is an embedded copy of $\mathfrak I$ with $\beta \times \{0\} \subset S,\ \beta \times \{1\} \subset R$, and $\mathfrak I \cap (S\cup R) = \partial \mathfrak I$.
Similarly, if $W$ is a compact surface with boundary and $\delta$ is a closed $1$--submanifold of $\partial W$, we denote by $\mathfrak J$ the product $W \times I$ with intervals $\{x\} \times I$ collapsed to points for $x\in \delta$.  A \emph{product region} in $M$ (called a \emph{blister} in \cite{Sa}) is an embedded copy of $\mathfrak J$ with $W \times \{0\} \subset S,\ W \times \{1\} \subset R$, and $\mathfrak J \cap (S\cup R) = \partial \mathfrak J\setminus \text{int} (\mathfrak J \cap \partial M)$. Note that $\delta$ is allowed to be empty, in which case the product region is really a product.

We say that two surfaces $S,R$ in a manifold $M$ have \emph{simplified intersection}, if they do not bound any product region. In particular, if a component $\dot{S}$ of $S$ is isotopic to a component $\dot{R}$ of $R$, then we must have $\dot{S}=\dot{R}$.

We say that $S$ and $R$ are \emph{almost transverse} if for each component $\dot{S}$ of $S$ and $\dot{R}$ of $R$ either $\dot{S}$ equals $\dot{R}$ or they intersect transversely. In particular, if $S$ equals $R$ then $S$ and $R$ are almost transverse.

We say that surfaces $S$ and $R$ are \emph{almost disjoint} if for intersecting components $\dot{S}$ of $S$ and $\dot{R}$ of $R$ we have $\dot{S}=\dot{R}$. In particular, $S$ is almost disjoint from itself.

Note that for a pair of surfaces $S,R$, the surface $R$ can be
always isotoped to $R'$ which is almost transverse to $S$ and has
simplified intersection with $S$. (This is not true if we wanted
to drop 'almost', consider the case where some components of $S$
and $R$ coincide. Actually, this also fails in the very special
case where $S=R$ and $M$ is a surface bundle over a circle, but we
will ignore that since then $\M$ is trivial.) Moreover, if $R_1,
R_2$ are almost disjoint, then they can be isotoped to almost
disjoint $R_1', R_2'$ which are both almost transverse to $S$ and
have simplified intersection with $S$ (again we cannot require
that $R_1', R_2'$ are disjoint, even if $R_1,R_2$ are).

\begin{rem}
In \cite{O} the definition of having simplified intersection consists of one more condition, which under standard hypotheses follows from the others. Namely, let $M$ be orientable, irreducible, $\partial$--irreducible and suppose that $S,R$ are orientable, incompressible and $\partial$--incompressible. If $S$ and $R$ are almost transverse and have simplified intersection, then there are no components of $S\cap R$ which are closed curves that are trivial in $S$ or $R$, or arcs that are $\partial$--parallel in $S$ or $R$.
\end{rem}

We now answer the opening question of the section.

\begin{prop}
\label{prop:simplified}
Let $S,R$ be spanning surfaces in $E$ representing $\sigma, \rho$ in $\MV$. If $S$ and $R$ are almost transverse and have simplified intersection, then they satisfy $$d(\sigma, \rho)=d(S,R).$$
\end{prop}

We deduce Proposition~\ref{prop:simplified} from the following
version of \cite[Proposition 4.8(2)]{Sa}, which we give without a
proof.

\begin{prop}
\label{prop:empty inters}
Let $M$ be (possibly non-compact) orientable, irreducible, and $\partial$--irreducible $3$--manifold. Let $W,N$ be (possibly non-compact) proper $3$--submanifolds of $M$ such that $\partial W, \partial N$ are incompressible and $\partial$--incompressible surfaces which are almost transverse with simplified intersection. If $N$ is isotopic to a submanifold $N'$ such that the interior of $N'$ is disjoint from $W$, then also the interior of $N$ is disjoint from $W$.
\end{prop}

In the setting described in Section~\ref{sec:distance}, this yields the following.

\begin{cor}
\label{cor:empty inters}
Let $W,N$ be proper $3$--submanifolds of $\widetilde{E}$ such that $\partial W, \partial N$ are unions of lifts of minimal genus spanning surfaces which are almost transverse with simplified intersection. If $N$ is isotopic to $N'$ such that the interior of $N'$ is disjoint from $W$, then also the interior of $N$ is disjoint from $W$.
\end{cor}

We will usually invoke Corollary~\ref{cor:empty inters} in the situation where $W=\overline{E}_j$ and
$N=\tau^i(\overline{E}^R)$ for some $j,i$, where $E_j$ and $E^R$ are as in Section~\ref{sec:distance}.

\medskip

We are now prepared for the following.

\begin{proof}[Proof of Proposition~\ref{prop:simplified}.]
Let $R$ and $S$ be almost transverse with simplified intersection. Let $R'$ be an element of $\rho=[R]$ for which the minimum of $d(S,R')$ is attained.
Then we have $d(\sigma, \rho)=d(S,R')=r'-m',$ where $E^{R'}, r', m'$ are as in Definition~\ref{def:distance} with $R$ replaced by $R'$. Then $E^{R'}$ is disjoint from all $E_{j}$ with $j\geq r'+1$ or $j\leq m'-1$. Let $E^R$ be the lift of $E\setminus R$ to $\widetilde{E}$ isotopic to $E^{R'}$.
Since $R$ has simplified intersection with $S$, its lifts have simplified intersection with the lifts of $S$. By Corollary~\ref{cor:empty inters}, $E^{R}$ is disjoint from all $E_{j}$ with $j\geq r'+1$ or $j\leq m'-1$. Then we have $r\leq r'$ and $m\geq m'$, which implies $d(S,R)\leq d(S,R')$, as desired.
\end{proof}

We conclude with recording the following lemma, whose proof we leave for the reader.

\begin{lemma}
\label{lem:triple}
Let $M$ be orientable, irreducible, $\partial$--irreducible and suppose that $S,R$ and $T$ are orientable, incompressible and $\partial$--incompressible surfaces properly embedded in $M$. Then $S,R$ and $T$ can be isotoped to be pairwise almost transverse and have pairwise simplified intersection.
\end{lemma}

\section{Projection map}
\label{sec:projection}
In this section we recall a construction of Kakimizu which we think of as a projection map and which will be our main tool. First, we need to fix a basepoint $\sigma\in \MV$.
The projection map $\pi_\sigma$ will map every $\rho\in \MV$ at distance $n>0$ from $\sigma$ to a vertex $\pi_\sigma(\rho)\in \MV$ adjacent to $\rho$ at distance $n-1$ from $\sigma$.

The existence of such projection map completes Kakimizu's proof of Proposition~\ref{prop:distances coincide}. It implies, in particular, that $\M$ is connected. In the present article we promote this method to prove contractibility of $\M$.

\medskip

We say that an oriented surface $T$ is obtained by a \emph{cut-and-paste operation} on $S$ and $R$ if it is a union of closures of oriented components of $S\setminus R,\ R\setminus S$ and common components of $S$ and $R$, with $\partial T\subset \partial S\cup \partial R$.

\begin{defin}
\label{def:projection}
Let $\sigma\neq\rho$ be vertices of $\M$. Put $n=d(\sigma, \rho)$.
For any fixed spanning surface $S\in \sigma$ we can choose $R\in \rho$ such that $S$ and $R$ are almost transverse with simplified intersection. In particular $S$ and $R$ have \emph{almost disjoint} boundaries, which means that the boundary components are disjoint or equal. By Proposition~\ref{prop:simplified} we have $d(S,R)=n$.

Recall the notation of Section~\ref{sec:distance} that $r$ is
largest such that the translate $E_r$ of $E_0$ intersects the lift
$E^R$ of $E\setminus R$ to $\widetilde{E}$. Denote
$\widetilde{R}=\overline{E}^R\cap \tau(\overline{E}^R)$. Let
$\widetilde{P}\subset S_r\cup \widetilde{R}$ denote the surface
obtained by a cut-and-paste operation on $S_r$ and
$\widetilde{R}$, which is the intersection of the boundaries of
$\overline{E}^R\setminus E_r$ and $\tau(\overline{E}^R)\cup
\overline{E}_r$. See Figure 2.

\begin{figure}[htbp]
\begin{center}
\includegraphics[width=5in]{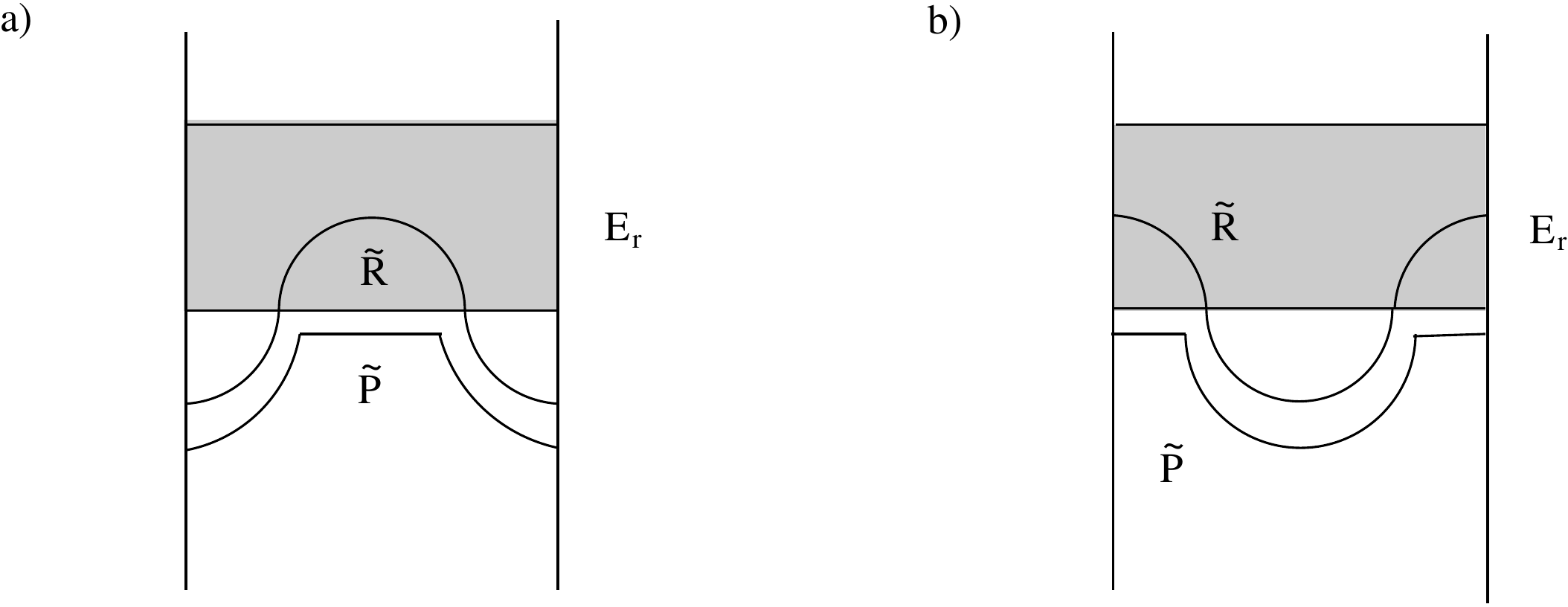}
\caption{Construction of $\widetilde P$}
\end{center}
\end{figure}

The surface $\widetilde{P}$ considered with the orientation
inherited from $\widetilde{R}$ and $S_r$ satisfies in homology
$\partial (E^R\cap E_r)=\widetilde{R}-\widetilde{P}$. Hence the
image $P$ of $\widetilde{P}$ under $p$ is in the homology class
$\alpha$. Moreover, $\widetilde{P}$ embeds under $p$ into $E$. In
order to justify that $P$ is a~spanning surface, it remains to
prove that its boundary $\partial P$ is not only homological but
also homotopic to $\gamma$. This follows from the fact that
$\partial P$  is homotopic to a~combination of curves in $\gamma$
and that by the hypothesis that $\gamma$ does not separate the
components of $\partial E$ no non-trivial combination of curves in
$\gamma$ is homological to zero.

Now a calculation as in case 1 of the proof of \cite[Theorem 2.1]{K} yields
that $P$ is a spanning surface of minimal genus. We define $$\pi_\sigma(\rho)=[P].$$ We prove that this
class is well-defined in Proposition~\ref{prop:well-defined}.
\end{defin}

As indicated at the beginning of this section, we have the following property, which justifies calling $\pi_\sigma$ the \emph{projection}.

\begin{rem}
\label{rem:projection basic}
The surface $P$ in Definition~\ref{def:projection} satisfies $d(R,P)=1$ and $d(S,P)=n-1$. Hence $\pi_\sigma(\rho)$ is adjacent to $\rho$ and satisfies $d(\sigma, \pi_\sigma(\rho))=n-1$.
\end{rem}

In the proof that the projection is well-defined we need the following result.
\begin{theorem}[{\cite[Theorem 3]{O}}]
\label{thm:Oertel}
Let $M$ be an orientable, irreducible, $\partial$--irreducible $3$--manifold. Let $S,R$ be orientable, incompressible, $\partial$--incompressible surfaces properly embedded in $M$. Assume that $S$ and $R$ are almost transverse with simplified intersection and that they are isotopic to $S',R',$ respectively, which are also almost transverse with simplified intersection. Suppose a cut-and-paste operation on $S$ and $R$ yields an orientable, incompressible and $\partial$--incompressible surface $P$. Then there is a corresponding cut-and-paste operation on $S',R'$ yielding a surface $P'$ isotopic to $P$.
\end{theorem}

\begin{prop}
\label{prop:well-defined}
The class $[P]$ in Definition~\ref{def:projection} does not depend on the choice of $S$ and $R$.
\end{prop}

\begin{proof}
We can fix $S\in \sigma$. Let $R,R'\in \rho$ be almost transverse to $S$ with simplified intersection.
Let $\widetilde{P}$ be obtained by a cut-and-paste operation on $\widetilde{R}$ and $S_r$ as in Definition~\ref{def:projection}. Let $E^{R'},\ \widetilde{R}'$ be the lifts of $E\setminus R',\ R'$ to $\widetilde{E}$ isotopic to $E^{R},\ \widetilde{R},$ respectively. By Corollary~\ref{cor:empty inters}, $r$ is the largest integer such that $E^{R'}$ intersects $E_r$. Let $\widetilde{P}'$ be the surface obtained from the cut-and-paste operation on $S_r$ and $\widetilde{R}'$ described in Definition~\ref{def:projection}, with $\widetilde{R}'$ in place of $\widetilde{R}$.

By Theorem~\ref{thm:Oertel} there is a surface $\widetilde{P}''$,
obtained by a cut-and-paste operation on $\widetilde{R}'$ and
$S_r$, which is isotopic to $\widetilde{P}$. The correspondence in
Theorem~\ref{thm:Oertel} (arising from the proof) is such that in
fact we have $\widetilde{P}''=\widetilde{P}'$, as desired.
\end{proof}

\section{Ordering the vertices}
\label{sec:order}

In this section we describe a natural way of ordering the vertices of the complex $\M$. One can check
that for special arborescent links this order coincides with the order described 
in \cite[Lemma 3.7]{Sa} (for appropriate $\sigma$).

\medskip

We begin with the following, which describes a possible position
of a pair of adjacent vertices $\rho,\rho'\in \MV$ with respect to
a vertex $\sigma\in \MV$. Note that $\rho$ and $\rho'$ may be at
the same or different distance from $\sigma$. We may choose almost
disjoint $R\in \rho, R'\in \rho'$ such that $R$ and $R'$ are
almost transverse to a~fixed $S\in \sigma$ and have simplified
intersection with $S$. Moreover, we can assume that $R$ and $R'$
have also simplified intersection (this does not follow
automatically from almost disjointness). By
Proposition~\ref{prop:simplified} we then have $d(R,R')=1$. As
usual $E^{R'}$ denotes a lift of $E\setminus R'$ to
$\widetilde{E}$ and $r'$ is largest such that $E_{r'}$ intersects
$E^{R'}$. Let $E^{R}$ be the lift of $E\setminus R$ contained in
$\overline{E}^{R'}\cup \tau^{-1}(\overline{E}^{R'})$.

\begin{defin}
\label{def:order} If $E^{R}$ intersects $E_{r'}$, then we write
$$\rho<_\sigma \rho'.$$ See Figure 3. We write $\rho\leq_\sigma \rho'$ if
$\rho<_\sigma \rho'$ or $\rho=\rho'$.
\end{defin}

\begin{figure}[htbp]
\begin{center}
\includegraphics[width=2in]{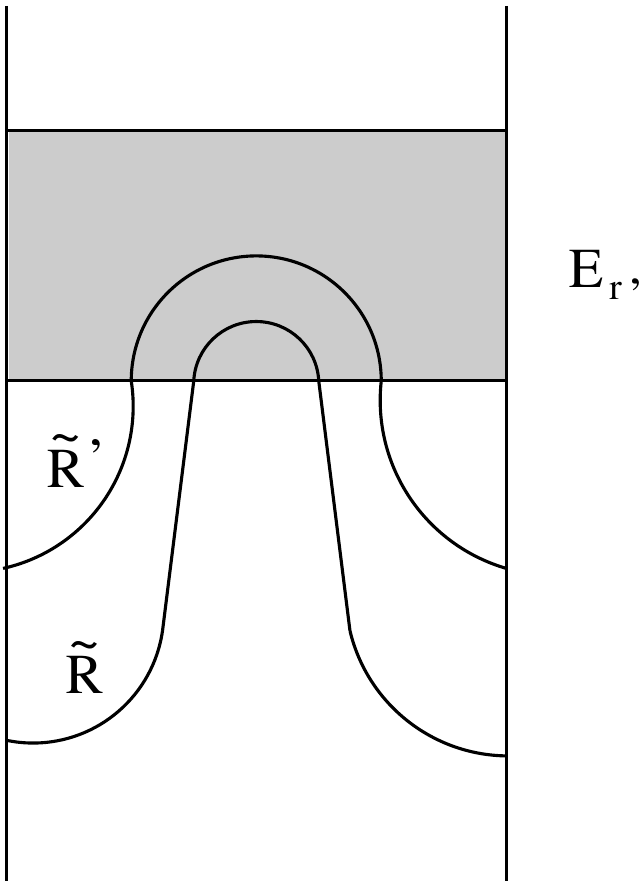}
\caption{Relation $[R] <_{\sigma} [R']$}
\end{center}
\end{figure}

\begin{rem}
Definition~\ref{def:order} does not depend on the choices of $R$
and $R'$. Indeed, by Corollary~\ref{cor:empty inters} the isotopy
class of $E^{R'}$ does not depend on the choice of $R'\in \rho'$.
Hence also the isotopy class of $E^{R}$ is well-defined. Again by
Corollary~\ref{cor:empty inters} the property that $E^{R}$
intersects $E_{r'}$ is invariant.
\end{rem}

We prove that adjacent vertices are always related by $<_\sigma$.

\begin{lemma}
\label{lem:neighbors related}
Let $\rho\neq \rho'$ be adjacent vertices of $\M$ and consider any $\sigma\in \MV$.
Then we have $\rho'<_\sigma \rho$ or $\rho<_\sigma \rho'$.
\end{lemma}

Later, in Lemma~\ref{lem:no cycles}, we will show that in fact $\rho'<_\sigma \rho$ and $\rho<_\sigma \rho'$ cannot happen simultaneously, which justifies using the notation $<_\sigma$.

\begin{proof}
Assume we do not have $\rho<_\sigma \rho'$, i.e.\ $E^R$ is disjoint from $E_{r'}$. If we now interchange $\rho$ with $\rho'$, then $E_r=E_{r'-1}$ takes on the role of $E_{r'}$ and $\tau^{-1}(E^{R'})$ takes on the role of $E^R$. Since $\tau^{-1}(E^{R'})$ intersects $E_{r'-1}$, we have $\rho'<_\sigma \rho$.
\end{proof}

In the following configuration we can determine the direction of the relation $<_\sigma$.

\begin{lemma}
\label{lem:diff distance}
If in Definition~\ref{def:order} the vertex $\rho$ is farther from $\sigma$ then $\rho'$, then we have $\rho<_\sigma \rho'$.
\end{lemma}
\begin{proof}
Since $E^{R}$ is contained in $\overline{E}^{R'}\cup \tau^{-1}(\overline{E}^{R'})$, it may intersect only $E_k$ with $m'-1\leq k \leq r'$. By Proposition~\ref{prop:simplified} we have $d(S, R)=d(S,R')+1$, so $E^{R}$ must intersect all those $E_k$. In particular it intersects $E_{r'}$, as desired.
\end{proof}

We now prove that, in particular, $\rho'\leq_\sigma \rho$ and $\rho\leq_\sigma \rho'$ implies $\rho=\rho'$.

\begin{lemma}
\label{lem:no cycles} There are no $\rho^1,\ldots, \rho^k,$ for
$k\geq 2,$ satisfying $$\rho^1<_\sigma \rho^2<_\sigma\ldots
<_\sigma \rho^k<_\sigma \rho^1.$$
\end{lemma}

Before we provide the proof, we record the following immediate
consequence. Note that in general the relation $<_\sigma$ is not
transitive, because $\rho<_\sigma \rho'$ and $\rho'<_\sigma
\rho''$ do not imply that $\rho$ and $\rho''$ are adjacent.

\begin{cor}
\label{lem:order}
The relation $<_\sigma$ extends to a linear order on $\MV$.
\end{cor}

\begin{proof}[Proof of Lemma~\ref{lem:no cycles}]
Since consecutive $\rho^i$ are adjacent, we can inductively choose
$R^k\in \rho^k,R^{k-1}\in \rho^{k-1},\ldots,R^1\in \rho^1$
satisfying the following. First, each $R^i$ is almost transverse
to $S$ with simplified intersection. Second, for $i<k$ the surface
$R^i$ is almost disjoint with $R^{i+1}$ and they have simplified
intersection. Let $r$ be largest such that $E_r$ intersects a lift
$E^{R^k}$ of $E\setminus R^k$. For $i<k$ define inductively
$E^{R^i}$ to be the lift of $E\setminus R^i$ contained in
$\overline{E}^{R^{i+1}}\cup \tau^{-1}(\overline{E}^{R^{i+1}})$. In
view of $\rho^1<_\sigma \rho^2<_\sigma\ldots <_\sigma \rho^k,$ all
$E^{R^i}$ intersect $E_r$.

Finally, let $R^*\in \rho^k$ be almost transverse to $S$ with simplified intersection and almost disjoint from $R^1$ with simplified intersection. Let $E^{R^*}$ be the lift of $E\setminus R^*$ contained in $\overline{E}^{R^{1}}\cup \tau^{-1}(\overline{E}^{R^{1}})$. In view of $\rho^k<_\sigma \rho^1$, $E^{R^*}$ intersects $E_r$. By Corollary~\ref{cor:empty inters}, $E^{R^*}$  and $E^{R^1}$ lie in the same isotopy class. Then the surfaces $\overline{E}^{R^*}\cap \tau(\overline{E}^{R^*})$ and $\overline{E}^{R^1}\cap \tau(\overline{E}^{R^1})$ are almost disjoint and bound a product containing all $\overline{E}^{R^i}\cap \tau(\overline{E}^{R^i})$. Hence all $\rho_i$ coincide, contradiction.
\end{proof}

\section{Properties of the projection map}
\label{sec:properties}

In this section we collect the properties of the projection map which will be later used to prove the theorems from the Introduction.

\medskip

The following property of the projection map $\pi_\sigma$ is the key to our proof of Theorem~\ref{thm:Contractibility}.

\begin{lemma}
\label{lem:same layer} Let $\rho$ and $\rho'$ be adjacent vertices
of $\M$ such that $\rho$ is different from some $\sigma\in \MV$.
Assume $\rho\leq_\sigma \rho'$. Then we have $\rho' \leq_\sigma
\pi_\sigma(\rho)$. In particular $\pi_\sigma(\rho)$ and $\rho'$
are equal or adjacent.
\end{lemma}

\begin{proof}
Let $S,R,R',E^{R'},E^{R}$ be as in Definition~\ref{def:order} and let
$\widetilde{P}$ be as in Definition~\ref{def:projection}. Let $E^P$ be that lift of $E\setminus P$ which contains $E^R\setminus E_r$.

Then $E^{R'}$ is contained in $\overline{E}^{P}\cup
\tau(\overline{E}^{P})$. In particular $\pi_\sigma(\rho)$ and
$\rho'$ are equal or adjacent. There is an isotopy $i$ of $P$ such
that $i(P)$ is almost transverse to $S$ with simplified
intersection and almost disjoint with $R'$ with simplified
intersection. Since $E^P$ is disjoint from $E_r$, by
Corollary~\ref{cor:empty inters} so is the lift of $E\setminus
i(P)$ in the isotopy class of $E^P$. Hence we do not have
$\pi_\sigma(\rho)<_\sigma \rho'$. By Lemma~\ref{lem:neighbors
related} we then have $\rho'\leq_\sigma\pi_\sigma(\rho)$, as
desired. See Figure 4.
\end{proof}

\begin{figure}[htbp]
\begin{center}
\includegraphics[width=3in]{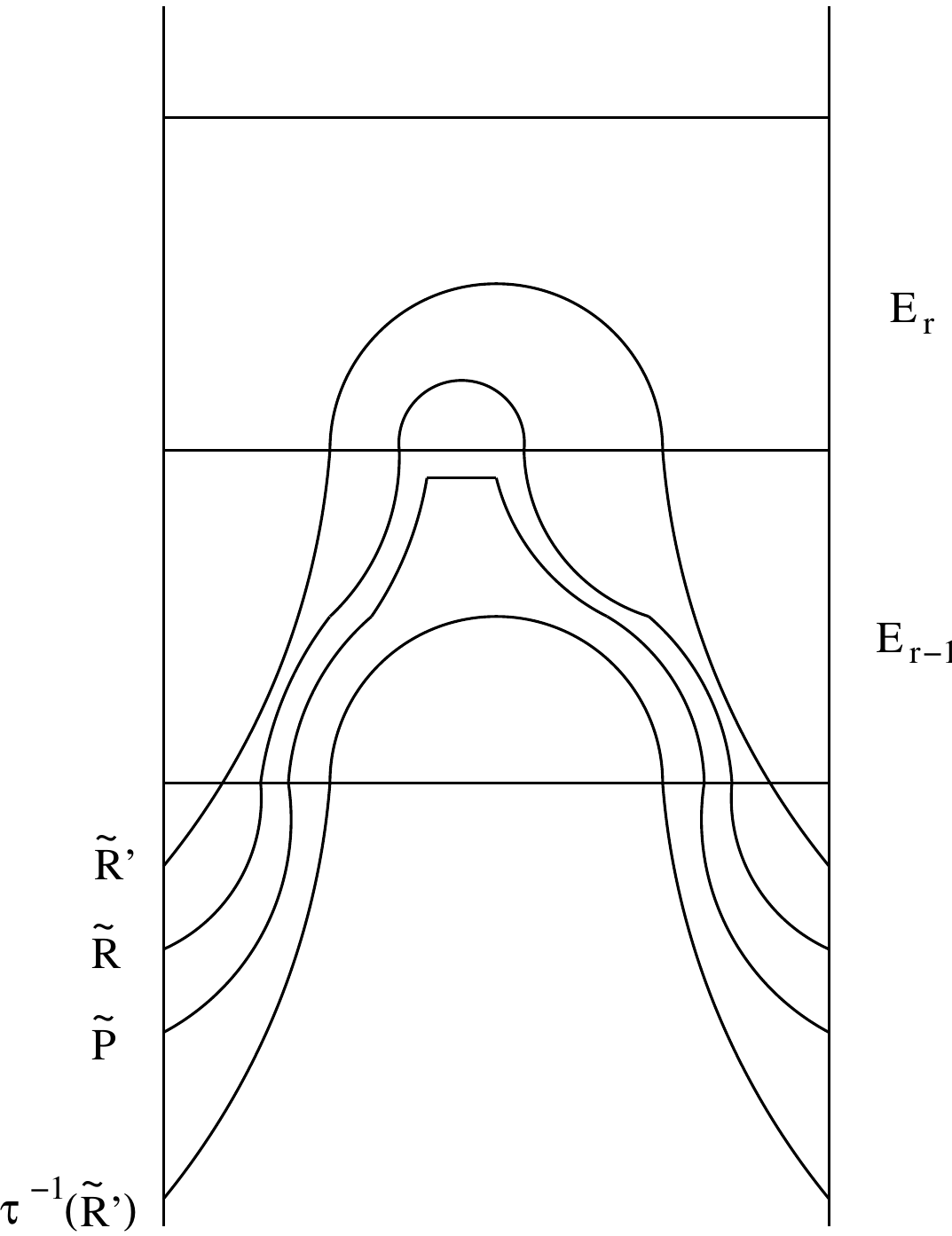}
\caption{Configuration from Lemma \ref{lem:same layer}}
\end{center}
\end{figure}

A double application of Lemma~\ref{lem:same layer} yields the following.

\begin{cor}
\label{lem:two projections}
Let $\rho$ and $\rho'$ be adjacent vertices of $\M$ different from some $\sigma\in \MV$. Assume $\rho\leq_\sigma \rho'$. Then we have $\pi_\sigma(\rho)\leq_\sigma \pi_\sigma(\rho')$.
\end{cor}

The following two results will be only used in the proof of Theorem~\ref{thm:fixed point} in Section~\ref{sec:fixed}. They are inspired by \cite{Pol}. In particular the proof of our Lemma~\ref{lem:bound in layer} resembles the proof of \cite[Lemma~3.9]{Pol}.

\begin{figure}[htbp]
\begin{center}
\includegraphics[width=3in]{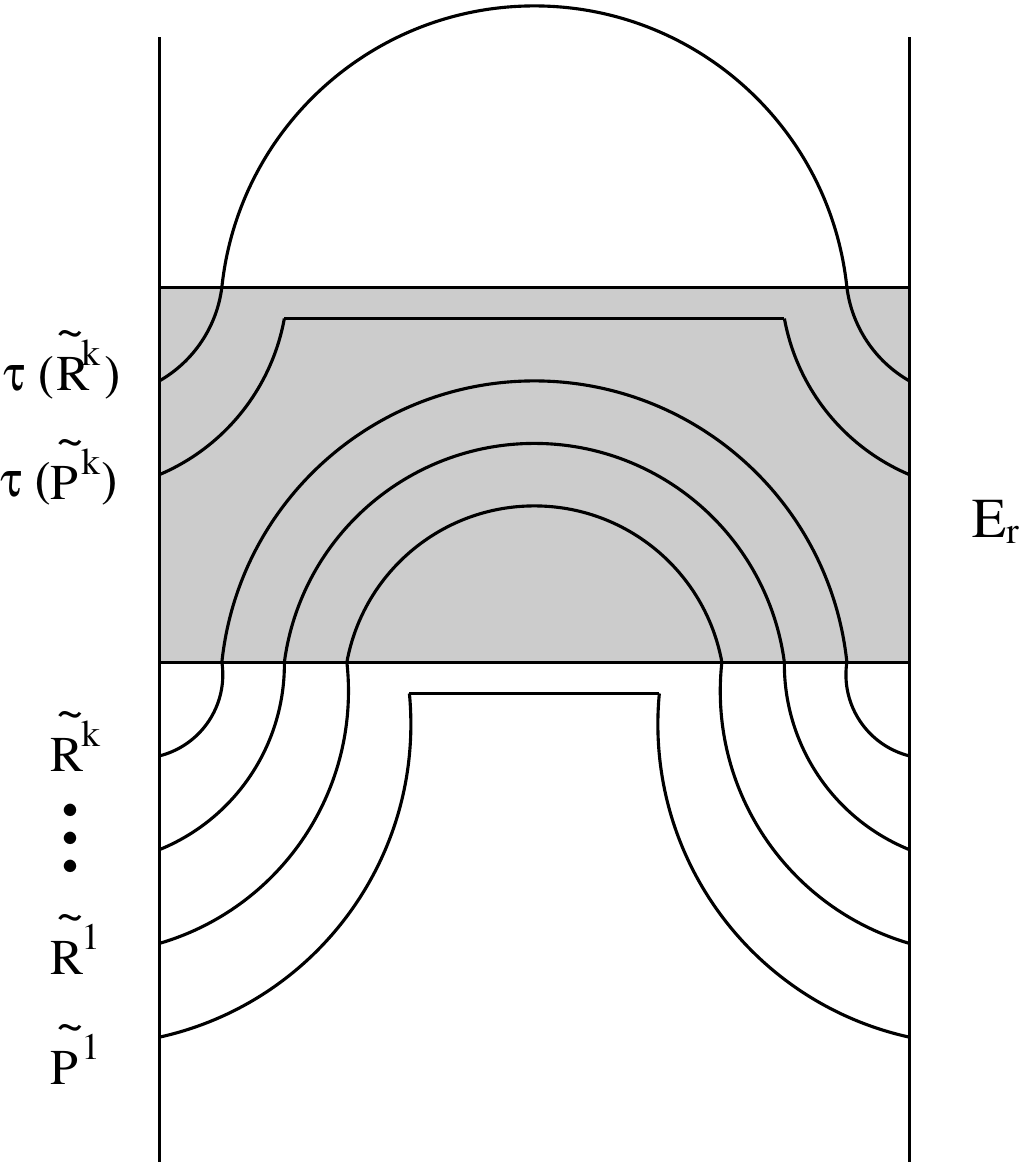}
\caption{Configuration from Lemma \ref{lem:same projection}}
\end{center}
\end{figure}

\begin{lemma}
\label{lem:same projection}
Assume that there are vertices $\rho^1,\ldots, \rho^k$ at the same non-zero distance from $\sigma\in\MV$ satisfying $$\rho^1<_\sigma \rho^2<_\sigma\ldots <_\sigma \rho^k,\text{ and } \pi_\sigma(\rho^1)=\pi_\sigma(\rho^k).$$
Then all $\pi_\sigma(\rho^i)$ are equal and all $\rho^i$ are adjacent.
\end{lemma}
\begin{proof}
The fact that all $\pi_\sigma(\rho^i)$ are equal follows
immediately from Corollary~\ref{lem:two projections} and
Lemma~\ref{lem:no cycles}. To show that all $\rho^i$ are adjacent
it is enough to give an argument that $\rho^1$ and $\rho^k$ are
adjacent (for other pairs of $\rho^i$ we pass to a subsequence).

First we choose $E^{R^k},\ldots,E^{R^1}$ in the same way as in the
proof of Lemma~\ref{lem:no cycles}. Let $\widetilde{P}^1,
\widetilde{P}^k$ be obtained as in
Definition~\ref{def:projection}. Then $\tau(\widetilde{P}^k)$ is
disjoint from $E^{R^k}$ and in the same isotopy class as
$\tau(\widetilde{P}^1)$. See Figure 5. Hence $\widetilde{R}^1$ and
$\widetilde{R}^k$ are isotopic to almost disjoint surfaces
$i(\widetilde{R}^1)$ and $i(\widetilde{R}^k)$ contained in the
closure of the lift of $E\setminus P$ bounded by $\widetilde{P}^1$
and $\tau(\widetilde{P}^1)$. Then we have
$$d\big(p(i(\widetilde{R}^1)),p(i(\widetilde{R}^k))\big)=1.$$
\end{proof}

Recall that by \cite[Proposition 4.9(1)]{Sa} all simplices of $\M$ can be realised by sets of disjoint spanning surfaces. Hence by Kneser's theorem, there is a bound on the dimension of simplices in $\M$. We promote this to the following.

\begin{lemma}
\label{lem:bound in layer}
For any $n>0$ there is a constant $l_n$ satisfying the following. Let $\sigma$ be any vertex of $\M$ and let $\rho^1,\ldots, \rho^l$ be at distance $n$ from $\sigma$ satisfying $$\rho^1<_\sigma \rho^2<_\sigma\ldots <_\sigma \rho^l.$$
Then we have $l\leq l_n$.
\end{lemma}
\begin{proof}
Let $L$ be a bound on the dimension of simplices in $\M$. We prove
by induction that it suffices to put $l_n=(L+1)^n$. For $n=1$ this
follows directly from Lemma~\ref{lem:same projection}. Assume we
have verified this for some $n\geq 1$.

Let now $\rho^1,\ldots, \rho^l$ be at distance $n+1$ from $\sigma$
satisfying $\rho^1<_\sigma \rho^2<_\sigma\ldots <_\sigma \rho^l$.
Put $i_0=0$. For $k\geq 1$ define inductively $i_k$ to be maximal
satisfying $\pi_\sigma(\rho^{i_k})=\pi_\sigma(\rho^{i_{k-1}+1})$,
until some $i_m$ equals $l$. By Lemma~\ref{lem:same projection}
for all $1\leq k\leq m$ we have $i_k-i_{k-1}\leq L+1$. Summing up
this implies $l\leq m(L+1)$.

It remains to bound $m$. By Corollary~\ref{lem:two projections} for all $1\leq k<m$ we have $\pi_\sigma(\rho^{i_k})<_\sigma \pi_\sigma(\rho^{i_k+1})$. This gives rise to
$$\pi_\sigma(\rho^{i_1})<_\sigma \pi_\sigma(\rho^{i_2})<_\sigma\ldots<_\sigma \pi_\sigma(\rho^{i_m}).$$
By Remark~\ref{rem:projection basic}, all $\pi_\sigma(\rho^{i_k})$ are at distance $n$ from $\sigma$.
By induction hypothesis we have $m\leq l_n$. Altogether, $l$ is bounded by $l_{n+1}=l_n(L+1)$, as desired.
\end{proof}

We will also need in Section~\ref{sec:fixed} the following technical result. Roughly speaking it says that projection paths do not exit balls containing their endpoints.

\begin{lemma}
\label{lem:diam2}
For $\sigma\neq\rho, \sigma'\in \MV$ with $d(\sigma', \rho)\leq d$ and $d(\sigma', \sigma)\leq d$ we have $d(\sigma', \pi_\sigma(\rho))\leq d$.
\end{lemma}
\begin{proof}
Choose $S\in\sigma, R\in \rho, S'\in \sigma'$ which are pairwise
almost transverse with simplified intersection (see
Lemma~\ref{lem:triple}). Let $r, \widetilde{P},P$ be as in
Definition~\ref{def:projection}. Let $E^P$ be the lift of
$E\setminus P$ bounded by $\widetilde{P}$ and
$\tau^{-1}(\widetilde{P})$. Choose a lift $E'_0$ of $E\setminus
S'$ to $\widetilde{E}$ and denote $E'_k=\tau^k(E'_0)$.

\begin{figure}[htbp]
\begin{center}
\includegraphics[width=3in]{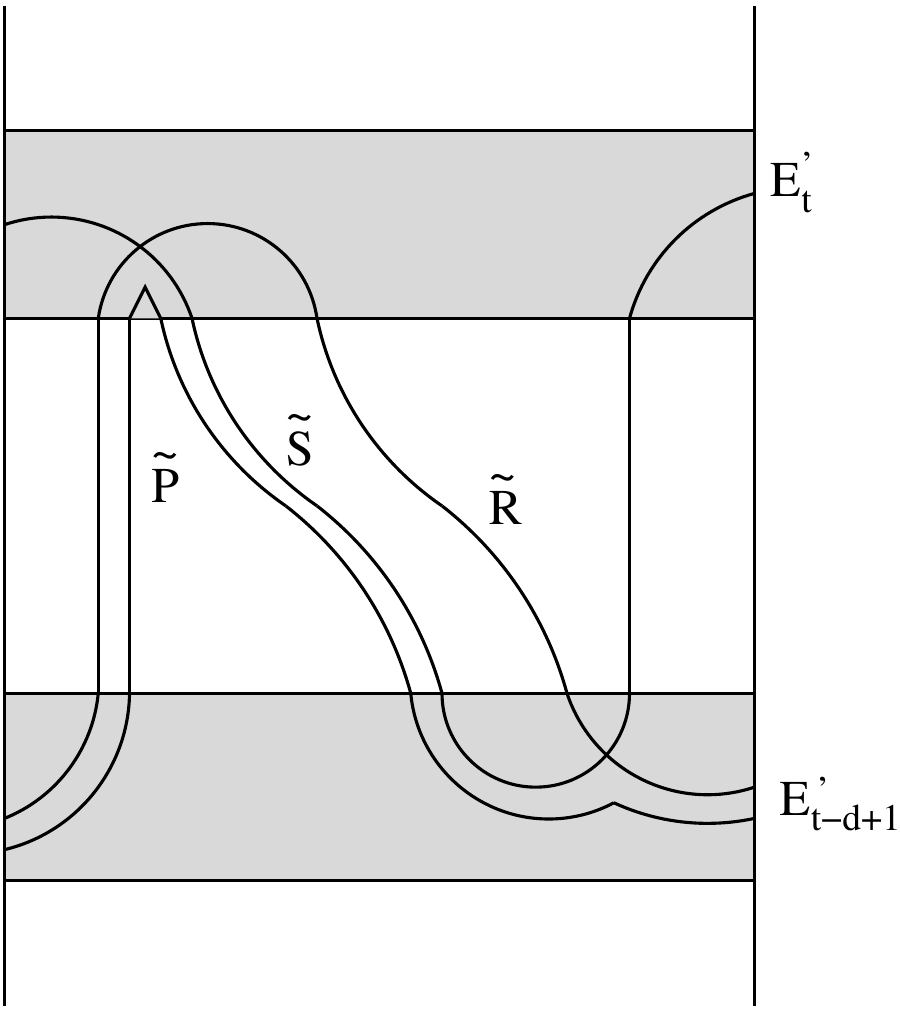}
\caption{Configuration from Lemma \ref{lem:diam2}  (here
$\widetilde S = \overline E_r \cap \overline E_{r+1}$)}
\end{center}
\end{figure}

Let $t$ be largest such that $E'_{t}$ intersects $E^R\cap E_{r-1}$
(which is non-empty). Note (see Figure 6) that $E^P$ is contained
in the union of
$$E^R\cap\big(\bigcup_{k\leq r-1}\overline{E}_{k}\big) \text{ and } E_{r-1}\cap \big(\bigcup_{i\leq
0}\tau^i(\overline{E}^R)\big).$$ In particular, $E^P$ is contained
in the intersection of $E^R\cup E_{r-1}$ with $\overline{E}'_k$'s
satisfying $k\leq t$. Since we have $d(S', R)\leq d$ and $d(S',
S)\leq d$, these $k$ must satisfy $t-k\leq d$, as desired.
\end{proof}

We conclude with another technical lemma which will be used only
in Section~\ref{sec:eg}. Roughly speaking, it describes how does
the projection $\pi_{\sigma'}$ look from the point of view of a
vertex $\sigma$ adjacent to $\sigma'$.

\begin{lemma}
\label{lem:change of basis}
Let $\sigma, \sigma'\in \MV$ be adjacent.
Let $\rho, \rho'\in \MV$ be also adjacent satisfying $\rho'<_{\sigma'} \rho$ and $\rho<_\sigma \rho'$. If $\sigma'\neq \rho'$, then we have
\begin{enumerate}[(i)]
\item
$\rho\leq_\sigma \pi_{\sigma'}(\rho')$,
\item
if $\sigma\neq \rho'$, then $d(\sigma,\pi_{\sigma'}(\rho'))\leq d(\sigma,\rho')$.
\end{enumerate}
\end{lemma}

See Figure 7 for an illustration.

\begin{figure}[htbp]
\begin{center}
\includegraphics[width=2.5in]{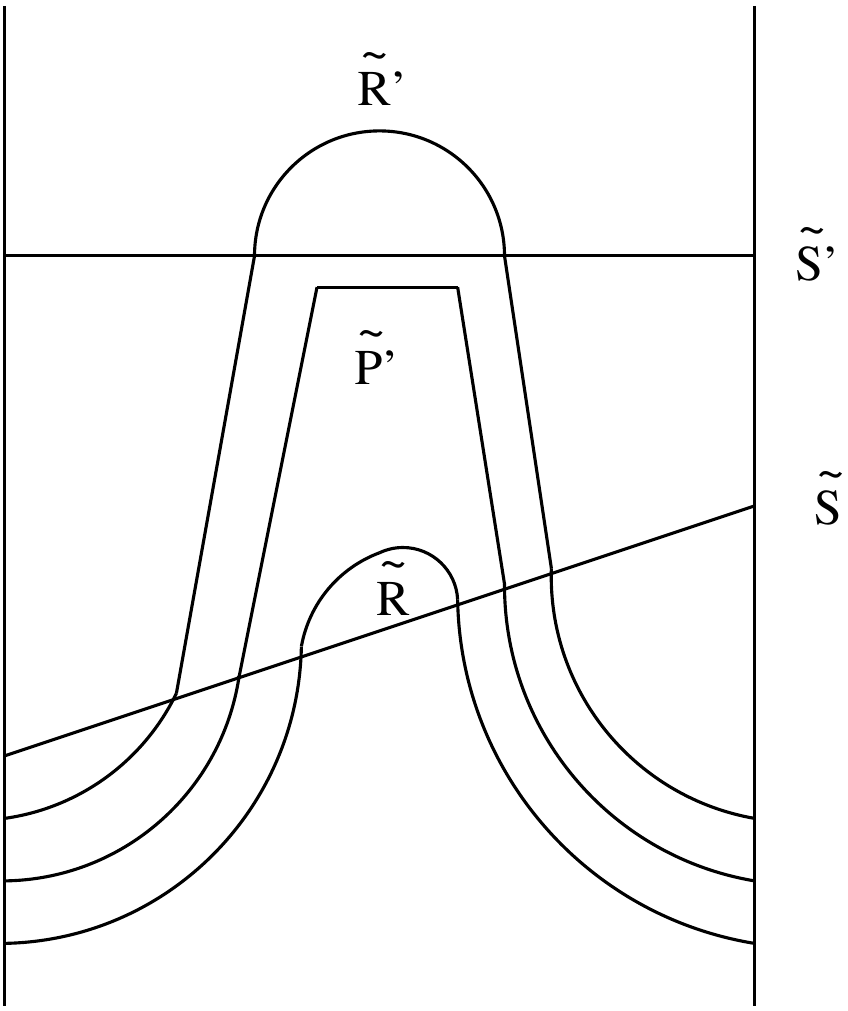}
\caption{Configuration from Lemma \ref{lem:change of basis} (here
$\widetilde S = \overline E_{r'} \cap \overline E_{r'+1}$ and
$\widetilde S' = \overline E'_{r'} \cap \overline E'_{r'+1}$)}
\end{center}
\end{figure}

\begin{proof}
Let $S\in \sigma, S'\in \sigma', R\in \rho, R'\in\rho'$ be
pairwise almost transverse with simplified intersection (this is
easily achieved by viewing $S\cup S'$ and $R\cup R'$ as a pair of
surfaces). Let $E_0'$ be the lift of $E\setminus S'$ contained in
$\overline{E}_0\cup \overline{E}_1$ (for some lift $E_0$ of
$E\setminus S$). Let $r'$ be largest such that
$E'_{r'}=\tau^{r'}(E'_0)$ intersects a lift $E^{R'}$ of
$E\setminus R'$. Let $E^R$ be the lift of $E\setminus R$ contained
in $\overline{E}^{R'}\cup \tau^{-1}(\overline{E}^{R'})$.

The hypotheses $\rho'<_{\sigma'} \rho$ and $\rho<_\sigma \rho'$
guarantee that $E^R$ is disjoint from $E'_{r'}$ but intersects
$E_{r'}$. Let $P'=p(\widetilde{P}')\in \pi_{\sigma'}(\rho')$ be
obtained as in Definition~\ref{def:projection} and let $E^{P'}$ be
the lift of $E\setminus {P'}$ bounded by $\widetilde{P}'$ and
$\tau^{-1}(\widetilde{P}')$. Since $E^R$ is disjoint from
$E'_{r'}$, the surface $\widetilde{P}'$ is contained in
$\tau(\overline{E}^R)$. (In particular $\rho$ and
$\pi_{\sigma'}(\rho')$ are equal or adjacent.)

There is an isotopy $i$ of $P'$ such that $i(P')$ is almost
transverse to $S$ with simplified intersection and almost disjoint
from $R$ with simplified intersection. Since $E^{P'}$ is disjoint
from $E_{r'+1}$, by Corollary~\ref{cor:empty inters} so is the
lift of $E\setminus i(P')$ in the isotopy class of $E^{P'}$.
Moreover, this lift contains $\widetilde{R}$ which intersects
$E_{r'}$. This implies assertion (i).

\medskip

Assertion (ii) is trivial since $E^{P'}$ intersects exactly the same $E_k$ as $E^{R'}$.
\end{proof}

\section{Contractibility}
\label{sec:contractibility}
In this section we prove Theorem~\ref{thm:Contractibility}. By Whitehead's theorem it suffices to prove that all finite subcomplexes of $\M$ are contained in contractible subcomplexes of $\M$.

We say that a flag subcomplex $X\subset \M$ is
\emph{$\sigma$--convex}, for $\sigma\in X^{(0)}$, if for any
$\rho\neq \sigma\in X^{(0)}$ we have $\pi_\sigma(\rho)\in
X^{(0)}$. By Remark~\ref{rem:projection basic} each finite
subcomplex of $\M$ is contained in a finite
\emph{$\sigma$--convex} subcomplex of $\M$ for any (hence some)
$\sigma$. Hence in order to prove
Theorem~\ref{thm:Contractibility}, it suffices to establish that
finite \emph{$\sigma$--convex} subcomplexes of $\M$ are
contractible. In fact, we have even a stronger property than
contractibility.

\begin{defin}
\label{def:dism}
A finite graph is \emph{dismantlable} if its vertices can be linearly ordered $x_0,\ldots, x_m$ so that for each $i\neq m$ there is $j>i$ satisfying
\begin{enumerate}[(i)]
\item
the vertex $x_j$ is adjacent to $x_i$,
\item
for any $x_k$ adjacent to $x_i$ with $k>i$, the vertex $x_j$ is adjacent or equal to $x_k$.
\end{enumerate}
\end{defin}

It is well known that finite flag complexes whose $1$--skeleta are dismantlable are contractible (see e.g.\ \cite{CO}). We just indicate that one obtains a homotopy retraction onto $x_m$ by successively retracting $x_i$ to $x_j$, where $j$ is as in Definition~\ref{def:dism}. In view of this in order to prove Theorem~\ref{thm:Contractibility} it remains to prove the following.

\begin{theorem}
\label{thm:dismantable}
Finite $\sigma$--convex subcomplexes of $\M$ have dismantlable $1$--skeleta.
\end{theorem}
\begin{proof}
We order all the vertices by extending the relation $<_\sigma$, which is possible by Corollary~\ref{lem:order}. By Lemma~\ref{lem:same layer} for all $\rho\neq \sigma$ we have $\rho<_\sigma\pi_\sigma(\rho)$, hence $\sigma$ is largest in this order.

For any non-largest $x_i$ we put $x_j=\pi_\sigma(x_i)$. As discussed above we have $x_i<_\sigma x_j$, which implies $j>i$ and condition (i) in Definition~\ref{def:dism}.

It remains to verify condition (ii). Let $x_k$ be adjacent to $x_i$ with $k>i$.
By Lemma~\ref{lem:neighbors related} we have $x_i<_\sigma x_k$ or $x_k<_\sigma x_i$. Since $k>i$ we must have $x_i<_\sigma x_k$. Then $x_j$ and $x_k$ are adjacent or equal by Lemma~\ref{lem:same layer}.
\end{proof}

\section{Fixed-point theorem}
\label{sec:fixed}
In this section we prove Theorem~\ref{thm:fixed point}. Key notions will be the following.

\begin{defin}
\label{def:semi}
A flag subcomplex $X$ of $\M$ is \emph{convex} if for all $\sigma\neq\rho\in X^{(0)}$ the vertex $\pi_\sigma(\rho)$ lies in $X^{(0)}$.

For a~vertex $v$ of $\M$, let $N(v)$ denote the union of $v$ with the set of all vertices adjacent to $v$. 
For a~subcomplex $X$ of $\M$ we put $N_X(v)=N(v)\cap X^{(0)}$.

A flag subcomplex $X$ of $\M$ is \emph{semi-convex} if for all $\sigma\neq\rho\in X^{(0)}$ 
there exists a vertex $\pi\in X^{(0)}$ 
satisfying
$$N_X(\pi_\sigma(\rho))\subset N_X(\pi),$$
and such that the distance between $\pi$ and $\sigma$ in the $1$--skeleton of $X$ equals 
$d(\pi_\sigma(\rho),\sigma)$. In particular, a convex subcomplex is also semi-convex.

The \emph{convex hull} of a subcomplex $X$ of $\M$ is the minimal
convex subcomplex of $\M$ containing $X$, i.e.\ it is the
intersection of all convex subcomplexes of $\M$ containing $X$.
\end{defin}

Note that semi-convex subcomplexes of $\M$ have $1$--skeleta
isometrically embedded in the $1$--skeleton of $\M$. Hence when we
discuss the distances in semi-convex subcomplexes we do not have
to specify if we consider the distance in the $1$--skeleton of the
subcomplex or of the whole $\M$. We also need the following
preliminary result which follows directly from
Lemma~\ref{lem:diam2}.

\begin{cor}
\label{lem:diam}
The convex hull of a subcomplex of diameter $d$ (in the $1$--skeleton of $\M$) has diameter $d$ as well.
\end{cor}

\begin{proof}[Proof of Theorem~\ref{thm:fixed point}]
Let $X\subset \MV$ be a finite orbit of the $G$--action on $\MV$. Denote by $\overline{X}$ the convex hull of $X$. By Corollary~\ref{lem:diam} $\overline{X}$ has finite diameter. Note that $\overline{X}$ is $G$--invariant. We consider now $G$--invariant non-empty semi-convex subcomplexes $Y$ of $\M$ of minimal diameter $d$. We want to show that $d$ equals $1$.

Otherwise, we also minimise the following value $l(Y)$. It is the maximum over $\sigma\in Y^{(0)}$ of $l$ admitting a sequence $\rho^1<_\sigma \rho^2<_\sigma\ldots <_\sigma \rho^l$ for some $\rho^1,\ldots, \rho^l$ at distance $d$ from $\sigma$. Note that $l(Y)$ is always finite by Lemma~\ref{lem:bound in layer}.

We say that a vertex $v$ of a subcomplex $Y$ of $\M$ is \emph{strongly dominated} (by $w$) in $Y$
if there is a vertex $w$ in $Y$ satisfying $N_Y(v)\subsetneq N_Y(w)$.

Let $Z$ denote the set of all the vertices $v\in Y^{(0)}$ strongly dominated in $Y$. 
Let $W$ be the subcomplex of $Y$ spanned by all the vertices in $Y^{(0)}\setminus Z$. 
Obviously $W$ is $G$--invariant. In order to obtain a contradiction it suffices to establish 
that $W$ is non-empty and semi-convex, and $l(W)<l(Y)$.

\medskip
We first prove $l(W)<l(Y)$. Consider any $\sigma\in W^{(0)}$ and a sequence $\rho^1<_\sigma \rho^2<_\sigma\ldots <_\sigma \rho^{l(Y)}$ of vertices at distance $d$ from $\sigma$. It suffices to show that $\rho^1$ belongs to $Z$.

By the definition of $l(Y)$ every $\rho\in Y^{(0)}$ at distance
$d$ from $\sigma$ adjacent to $\rho^1$ violates $\rho<_\sigma
\rho^1$. Then by Lemma~\ref{lem:neighbors related} we have
$\rho^1<_\sigma \rho$. By Lemma~\ref{lem:diff distance} the same
holds for all other $\rho\in Y^{(0)}$ adjacent to $\rho^1$. Hence
by Lemma~\ref{lem:same layer} all vertices in $Y$ adjacent to
$\rho^1$ are adjacent to or equal $\pi_\sigma(\rho^1)$. Note that
$\pi_\sigma(\rho^1)$ might note lie in $Y^{(0)}$ but since $Y$ is
semi-convex, there is $\pi\in Y^{(0)}$ at distance $d-1$ from $\sigma$ satisfying
$N_Y(\pi_\sigma(\rho^1))\subset N_Y(\pi)$. At
this point we have
$$N_Y(\rho^1)\subset N_Y(\pi).$$
Similarly, since $d\geq 2$, there is $\pi'\in Y^{(0)}$ 
at distance $d-2$ from $\sigma$
satisfying $N_Y(\pi_\sigma(\pi))\subset N_Y(\pi')$.
In particular, $\pi'$ is adjacent to $\pi$, but not to $\rho^1$. 
Hence we have $$N_Y(\rho^1)\subsetneq N_Y(\pi).$$
We conclude that $\rho^1$ is strongly dominated by $\pi$ in $Y$, which means that $\rho^1$ belongs to $Z$.

\medskip
We now prove that $W$ is non-empty. Pick a vertex $v\in Y^{(0)}$
with maximal $N_Y(v)$ (with respect to inclusion). Such
a vertex exists since otherwise we would have a simplex in $\M$ of
infinite dimension. Then $v$ is not strongly dominated in $Y$ by any
vertex and hence $v$ belongs to $W^{(0)}$.

\medskip
It remains to show that $W$ is semi-convex. Take
$\sigma\neq\rho\in W^{(0)}$. Since $Y$ is semi-convex, there is a~vertex $\pi$ of $Y^{(0)}$
at distance $d-1$ from $\sigma$ satisfying  $N(\pi_\sigma(\rho))\subset N(\pi)$. 
Let $\pi'$ be a~vertex of $Y^{(0)}$ with maximal possible $N_Y(\pi')$ containing $N_Y(\pi)$.
Such a vertex exists since $\M$ is finite-dimensional. Then $\pi'$ is not strongly dominated in $Y$, 
hence $\pi'$ belongs to $W^{(0)}$. Note that we also have $N_Y(\pi_\sigma(\rho))\subset N_Y(\pi')$.

Now we prove that $\pi'$ is at distance $d-1$ from $\sigma$ in $W^{(1)}$.
Let $\pi_0=\pi,\pi_1, \ldots, \pi_{d-1}=\sigma$ be a~
path in $Y^{(0)}$ from $\pi$ to $\sigma$. Put $\pi_0'=\pi',\pi_{d-1}'=\sigma$ and
for all $0<i< d-1$ let $\pi_i'$ be a~vertex of 
$Y^{(0)}$ with maximal possible $N_Y(\pi_i')$ containing $N_Y(\pi_i)$. Like before, all $\pi_i'$ belong to $W^{(0)}$.
Moreover, since $\pi_i$ is adjacent to $\pi_{i+1}$, also $\pi_i'$ is adjacent to $\pi_{i+1}$ and consequently $\pi_i'$
is adjacent to $\pi'_{i+1}$. Hence $\pi_i'$ form a~path and $\pi'$ is at distance $d-1$ from $\sigma$ in $W^{(1)}$.
Thus $W$ is semi-convex, as required.

\medskip
To summarise, assuming $d\geq 2$ we proved that $Y$ contains non-empty semi-convex 
$G$--invariant $W$ with $l(W)<l(Y)$ (where $l(W)=0$ means that the diameter of $W$ 
is less than $d$). This contradicts the choice of $Y$. In case $d=1$, $Y$ is the desired $G$--invariant simplex.
\end{proof}

Note that the proof would be easier if we knew that $\M$ is locally finite.

\section{Contractibility of fixed-point sets}
\label{sec:eg}

In this section we prove Theorem~\ref{thm:all contractible}. This is an elaboration on the proof from Section~\ref{sec:contractibility}.

Let $G$ be a subgroup of the mapping class group of $E$ fixing
$\alpha$ and the homotopy class of $\gamma$. Its fixed-point set
$\mathrm{Fix}_G(\M)$ has the following structure of a flag
simplicial complex $X$. Its vertices can be identified with the
set $V$ of minimal $G$--invariant simplices of $\M$. Its edges are
spanned on pairs vertices corresponding to simplices in $\M$
spanning a common simplex.

We assume that $X=\mathrm{Fix}_G(\M)$ is non-empty, i.e.\ there is a vertex $\Sigma\in V$ of $X$ (a simplex of $\M$) which is invariant under $G$. We need to prove that $X$ is contractible. The plan of the proof is the same as in Section~\ref{sec:contractibility}. We will define a mapping $\Pi_\Sigma$ from $V\setminus \{\Sigma\}$ to $V$ which will play the role of $\pi_\sigma$. We will observe that each finite subcomplex of $X$ lies
in a finite \emph{$\Sigma$--convex} subcomplex of $X$. The proof will then reduce to proving dismantlability of $\Sigma$--convex subcomplexes of $X$.

\begin{defin}
\label{def:big projection}
For $\Sigma\neq \Delta\in V$ we define $\Pi_\Sigma(\Delta)\in V$ in the following way.
We choose a vertex $\sigma$ of the simplex $\Sigma$. We consider $\delta\in \Delta$ which is minimal with respect to the order $<_\sigma$. We define $\Pi_\Sigma(\Delta)$ to be the $G$--orbit of $\pi_\sigma(\delta)$. We still need to check that this is an element of $V$, i.e.\ a simplex in $\M$. Note that since the relation $<_\sigma$ and the mapping $\pi_\sigma$ are $G$--equivariant, this definition does not depend on the choice of $\sigma$.
\end{defin}

\begin{lemma}
\label{lem:big projection well-defined}
$\Pi_\Sigma(\Delta)$ spans a simplex of $\M$. As a vertex of $X$ it is adjacent to $\Delta$.  Furthermore, for $\sigma\in \Sigma, \delta\in \Delta$ as in Definition~\ref{def:big projection} and all $\pi\in\Pi_\Sigma(\Delta)$ we have $$\delta\leq_\sigma \pi.$$
\end{lemma}
\begin{proof}
Let $\sigma\in \Sigma$ and $\delta\in \Delta$ be as in Definition~\ref{def:big projection}. By Lemma~\ref{lem:same layer}, for all $\delta'\in \Delta$ we have $\delta'\leq_\sigma \pi_\sigma(\delta)$. In particular, $\pi_\sigma(\delta)$ is adjacent or equal to all the vertices of $\Delta$.

Let now $\pi$ be any vertex of $\Pi_\Sigma(\Delta)$. By equivariance, $\pi$ is adjacent or equal to all the vertices of $\Delta$. Moreover, we have $\pi=\pi_{\sigma'}(\delta')$ for some $\sigma'\in \Sigma, \delta'\in \Delta$ satisfying $\delta'<_{\sigma'}\delta$. Now Lemma~\ref{lem:change of basis}(i) implies
$\delta\leq_\sigma\pi$.

Finally, by Lemma~\ref{lem:same layer}, $\pi_\sigma(\delta)$ and $\pi$ are adjacent or equal.
\end{proof}

We have the following analogue of Remark~\ref{rem:projection basic}, which in particular implies that $\Pi_\Sigma(\Delta)$ is different from $\Delta$.

\begin{lemma}
\label{lem:projection is nearer}
The sum of the distances between a vertex of $\Sigma$ and all the vertices of $\Pi_\Sigma(\Delta)$ is less than the corresponding sum for $\Sigma$ and $\Delta$.
\end{lemma}

Note that by equivariance the value in Lemma~\ref{lem:projection is nearer} does not depend on the choice of the vertex of $\Sigma$.

\begin{proof}
Fix $\sigma\in \Sigma$ and let $\delta\in \Delta$ be minimal with
respect to $<_\sigma$. By Remark~\ref{rem:projection basic} we
have $d(\sigma, \pi_\sigma(\delta))<d(\sigma, \delta)$. All other
vertices $\delta'\in\Delta$ are in correspondence with vertices
$\pi'\in \Pi_\Sigma(\Delta)$ of the form $\pi_{\sigma'}(\delta')$
for some $\sigma'\in \Sigma$. By Lemma~\ref{lem:change of
basis}(ii) we then have $d(\sigma, \pi')\leq d(\sigma, \delta')$.
Summing up the inequalities yields the lemma.
\end{proof}

We now introduce a definition analogous to the one in Section~\ref{sec:contractibility}.

\begin{defin}
A flag subcomplex $Y$ of $X$ is \emph{$\Sigma$--convex}, for $\Sigma\in Y^{(0)}$, if for any $\Delta\in Y^{(0)}\setminus \{\Sigma\}$ we have $\Pi_\Sigma(\Delta)\in Y^{(0)}$.
\end{defin}

Note that by Lemma~\ref{lem:projection is nearer} each finite subcomplex of $X$ is contained in a finite $\Sigma$--convex subcomplex of $X$. Hence in order to prove Theorem~\ref{thm:all contractible}, it remains to show the following.

\begin{theorem}
Let $Y$ be a finite \emph{$\Sigma$--convex} subcomplex of $X$. Then $Y^{(1)}$ is dismantlable.
\end{theorem}
\begin{proof}
We choose any $\sigma\in \Sigma$. By Corollary~\ref{lem:order} we can extend the relation $<_\sigma$ to a linear order on $\MV$. Let $x_0$ be the vertex of $Y^{(0)}$ containing the minimal possible vertex of $\M$ in this order. Let $x_1$ be one of the remaining vertices of $Y^{(0)}$ containing a minimal possible vertex of $\M$ etc. By Lemma~\ref{lem:big projection well-defined}, every $\Pi_\Sigma(\Delta)$ is larger than $\Delta$ in this order. In particular, $\Sigma$ is largest.

For any non-largest $x_i$ we put $x_j=\Pi_\Sigma(x_i)$. By Lemma~\ref{lem:big projection well-defined}  $j$ satisfies condition (i) in Definition~\ref{def:dism} and (as discussed above) we have $j>i$.

It remains to verify condition (ii). Let $x_k$ be adjacent to $x_i$ with $k>i$.
Let $\delta\in x_i$ be the minimal element with respect to $<_\sigma$. By the way we have ordered the $x$'s, for all $\delta'\in x_k$ we have $\delta<_\sigma \delta'$. From Lemma~\ref{lem:same layer} we get $\delta'\leq_\sigma \pi_\sigma(\delta)$, for all $\delta'\in x_k$. By equivariance, we get that $\delta'$ and $\pi$ are adjacent or equal for all $\delta'\in x_k$ and $\pi\in \Pi_\Sigma(x_i)=x_j$.
This means that $x_k$ and $x_j$ are adjacent or equal, as desired.
\end{proof}

\end{document}